\documentclass[10pt]{article}
\usepackage[utf8x]{inputenc}
\usepackage{amssymb,latexsym,amsmath,amsxtra,amsthm,amsfonts,amscd,slashed,mathrsfs,yhmath}
\usepackage{parskip}
\usepackage{enumerate}
\usepackage{color}
\usepackage{epsfig}
\usepackage{url}

\theoremstyle{plain}
\newtheorem{theorem}{Theorem}[section]
\newtheorem{corollary}[theorem]{Corollary}
\newtheorem{lemma}[theorem]{Lemma}

\theoremstyle{remark}

\theoremstyle{definition}

\makeatletter
\def\thm@space@setup{
  \thm@preskip=\parskip \thm@postskip=0pt
}
\makeatother

\newcommand{\paren}[1]{\left(#1\right)}
\newcommand{\norm}[1]{\Vert#1\Vert}

\newcommand{\abs}[1]{\left|#1\right|}

\def\R{\mathbb{R}}

\def\Z{\mathbb{Z}}
\def\J{\mathcal{J}}

\title{A Comparison of Bessel and Riesz Potentials.}
\date{}
\author{Ikemefuna Agbanusi}
\begin{document}

\maketitle

\begin{abstract}
How large is the Bessel potential, $G_{\alpha,\mu}f$, compared to the Riesz potential, $I_\alpha f$? 
 In this paper, we show that if $I_\alpha f\in L^p$ with $0<\alpha<1$ and $p>1$, then the following interpolation bound holds:
 \[\norm{G_{\alpha,\mu}f}_p\leq C(\omega(I_\alpha f,1/\mu)_p)^\alpha\cdot\norm{I_\alpha f}^{1-\alpha}_p.\]
Here $\omega(f,t)_p$ is the $L^p$ modulus of continuity. However, if $\alpha=p=1$, we obtain the ``$L\log L$" type result:
\[
\norm{G_{1,\mu} f}_{_1}\leq B\omega(I_1f,1/\mu)_1|\log\omega(I_1f,1/\mu)_1|.
\]
%whereas when $p>1$ and $I_1f\in L^p$, we have the equivalence
% \[\norm{G_{1,\mu} f}_p\approx \omega(I_1f,1/\mu)_p.\]
These and other estimates are obtained by studying the quotient of the two operators, $E_{\alpha,\mu}:=\frac{(-\Delta)^{\alpha/2}}{(\mu^2I-\Delta)^{\alpha/2}}$. This operator is of independent interest due to its connection to approximation theory. 

%Additionally, by analyzing its kernel, we show a  ``localization" result to the effect that $G_{\alpha,\mu}f(x_0)=\mathcal{O}(\mu^{-\frac\alpha2})$ if $I_\alpha f$ vanishes near $x_0$.
\end{abstract}

\section{Introduction}
Recall that if $\hat{f}(\xi)$ denotes the Fourier transform in $\R^d$, the kernel of the Bessel potential operator is defined by 
\[\widehat{G}_{\alpha,\mu}(\xi):=(\mu^2+|\xi|^2)^{-\frac\alpha2};\quad\alpha>0;\quad\mu>0,\] while the Riesz kernel is defined by
\[\widehat{I}_{\alpha}(\xi):=|\xi|^{-\alpha},\quad0<\alpha<d.\]
Both operators are indispensable in the theory of Sobolev spaces or, more broadly, spaces of functions with generalized derivatives. Their associated capacities are also vital for describing the fine structure of sets in various problems of analysis and PDE. 

The goal of this paper is to give yet another quantitative comparison of the two potentials with emphasis on the dependence on the parameter $\mu$. This touches on issues at the intersection of Fourier analysis, approximation theory and, of course, potential theory. 

The relationship between these two classical operators has been expressed in the literature in several ways. For instance, when $\mu=1$---which is the standard case---it is well known that the Riesz and Bessel kernels satisfy
\[
0<G_{\alpha,1}(x)\leq I_{\alpha}(x); \quad 0<\alpha<d.
\]
From this follows the convolution inequality $G_{\alpha,1}\star f(x)\leq I_\alpha\star f(x)$ which holds for $f\geq0$. This allows us to conclude that sets of Bessel capacity zero also have Reisz capacity zero as in {\sc Ziemer} \cite[p.\,67]{Ziemer:1989dp}. It also yields the pointwise comparison: if $f\geq0$ and $I_\alpha\star f(x_0)=0$ then
$G_{\alpha,1}\star f(x_0)=0$.
This property may not hold for functions which change sign, but we prove a result that says the Bessel potential cannot be large at points where the Riesz potential vanishes.

\begin{theorem}\label{thm:bessel-riesz-local}
If $I_\alpha f$ is in $L^p$ and vanishes near $x_0$, then $G_{\alpha,\mu} f(x_0) =\mathcal{O}(\mu^{-\frac\alpha2})$.  
\end{theorem}
The implicit constant depends on $p$ and the size of the neighborhood where $I_{\alpha}f$ vanishes.

The potentials have also been compared via norm estimates and we turn to describing these now. Already, the inequality $G_{\alpha,1}(x)\leq I_{\alpha}(x)$ yields $\norm{G_{\alpha,1}f}_p\leq \norm{I_{\alpha}f}_p$ if $f\geq0$ and $I_{\alpha}f\in L^p(\R^d)$. For $f\geq 0$, there is the deeper estimate
\begin{equation}\label{eqn:muc-whe-sch}
c\norm{I_{\alpha,\frac1\mu}f}_p\leq\norm{G_{\alpha,\mu}f}_p\leq C\norm{M_{\alpha,\frac1\mu}f}_p;\quad
\end{equation}
which holds for $ 1<p<\infty$, $0<\alpha<d$, constants $c,C>0$ and involves the \emph{truncated} Riesz kernel 
\[
I_{\alpha,\delta}f(x)=\int\limits_{|x-y|\leq\delta}f(y)I_\alpha(x-y)\,dy,
\]
as well as the fractional maximal operator
 \[
 M_{\alpha,\delta}f(x)=\sup_{\substack{Q\ni x\\ l(Q)\leq\delta}}\frac{1}{|Q|^{1-\frac\alpha d}}\int_Q|f(y)|dy.
\]
Here $Q$ denotes a cube with sides parallel to the coordinate planes and $l(Q)$ is its side length. For details see {\sc Adams--Hedberg} \cite[Theorem 3.6.2]{Adams:1999qc}, {\sc Schechter} \cite[Theorem 3.5]{Schecter:1989eu} and {\sc Muckenhoupt--Wheeden} \cite[Theorem 1]{Muckenhoupt:1974bd}.
 Our other main result is a version of \eqref{eqn:muc-whe-sch}. To state it, we need the $L^p$ modulus of continuity defined by
\begin{equation}\label{eqn:mod_cont}
\omega(f,t)_p:=\sup_{|h|\leq t}\norm{f(\cdot+h)-f(\cdot)}_p.
\end{equation}

\begin{theorem}\label{thm:bessel-riesz-modulus}
\leavevmode
\begin{enumerate}[(a)]
\item If $1<p<\infty$ and $I_1f\in L^p(\R^d)$, there is an $A>0$ such that
\[A^{-1}\omega(I_1f,1/\mu)_p\leq \norm{G_{1,\mu} f}_p\leq A\omega(I_1f,1/\mu)_p.\]
\item If $I_1f\in L^1(\R^d)$, then for some $B>0$
\[
\norm{G_{1,\mu} f}_{_1}\leq B\omega(I_1f,1/\mu)_1|\log\omega(I_1f,1/\mu)_1|.
\]
\item If $0<\alpha<1$,  $1\leq p<\infty$ and $I_\alpha f\in L^p(\R^d)$, there is a $C_\alpha>0$ such that
\[
\norm{G_{\alpha,\mu}f}_p\leq C_\alpha(\omega(I_\alpha f,1/\mu)_p)^\alpha\cdot\norm{I_\alpha f}^{1-\alpha}_p.
\]
\end{enumerate}
\end{theorem}

Theorem \ref{thm:bessel-riesz-modulus} thus gives a precise sense in which the Bessel potential fine--tunes the Riesz potential.

Here is a summary of the paper. The main idea is ultimately a simple one---to compare the two operators, we examine the ``Bessel--Riesz quotient":
\[
 E_{\alpha,\mu}:=\frac{(-\Delta)^{\alpha/2}}{(\mu^2I-\Delta)^{\alpha/2}}.
\]
The crux is that $G_{\alpha,\mu}f=E_{\alpha,\mu}I_\alpha f$, so estimates for $E_{\alpha,\mu}$ lead to estimates between two potentials. After establishing notation in \S2, we turn to Theorem \ref{thm:bessel-riesz-modulus} in \S3. Its proof exploits the approximation theoretic properties of the Bessel--Riesz quotient. Theorem \ref{thm:bessel-riesz-local} is proved in \S4 using Fourier analytic estimates on the symbol and kernel of $E_{\alpha,\mu}$.

It would be interesting to know if similar results hold for $\frac{\sqrt{L(x,D)}}{\sqrt{\mu^2I+L(x,D)}}$ where $L(x,D)$ is now a linear second order differential or pseudo-differential operator which need not be elliptic. We hope to tackle this in the future.

\section{Notation}
 Everything takes place in $d$--dimensional Euclidean space $\R^d$ and for $1\leq p\leq\infty$, $L^p=L^p(\R^d)$ are the usual Lebesgue spaces with norm denoted by $\norm{f}_p$. 

For a non-negative integer $k$, the \emph{Sobolev space} $W^k_p$ consists of $L^p$ functions having distributional derivatives up to order $k$ in $L^p$. In $W^k_p$ we use the norm $\norm{f}_{W^k_p}=\sum_{|\gamma|\leq k}\norm{D^\gamma f}_p$, and seminorm $|f|_{W^k_p}=\sum_{|\gamma|= k}\norm{D^\gamma f}_p$.

The direct and inverse Fourier transform of $f$ and $\hat{g}$ respectively defined as
\[\hat{f}(\xi)=\int e^{-ix\cdot\xi}f(x)\,dx;\qquad \check{g}(x)=(2\pi)^{-d}\int e^{ix\cdot\xi}\hat{g}(\xi)\,d\xi.\]
When convenient, we also use $\mathcal{F}_{x\to\xi}$ and $\mathcal{F}_{\xi\to x}$ for the direct and inverse transform. For suitable functions $a(\xi)$, we associate the multiplier operator
\[a(D)f(x)=(2\pi)^{-d}\int e^{ix\cdot\xi}a(\xi)\widehat{f}(\xi)\,d\xi.\]
We will often use the H\"ormander--Mikhlin multiplier theorem: if $|\partial^\gamma a(\xi)|\leq C_\gamma |\xi|^{-\gamma}$ for $|\gamma|\leq k$ with $k>d/2$, then $a(D)$ defines a bounded operator in $L^p$ for $1<p<\infty$. See  {\sc Stein} \cite[\S4.3.2]{Stein:1970yt} for details.

%The \emph{Bessel potential space} $\mathcal{L}_p^s$ is defined as $\mathcal{L}_p^s:=\{G_s\star f:f\in L^p\}$. Here $G_s(x):=G_{s,1}(x)$ is the ``standard" Bessel kernel.
% defined by \[G_s(x)=\mathcal{F}_{\xi\to x}\brac{(1+|\xi|^2)^{-\frac s2}}.\]
For $0<s<1$, $1\leq p<\infty$ and $1\leq q\leq\infty$, we define the \emph{Besov spaces} $B^s_{p,q}$ as those $f\in L^p$ for which the seminorm
\[|f|_{B^s_{p,q}}:=\paren{\int_0^1(t^{-s}\omega(f,t)_p)^q\frac{dt}{t}}^{1/q}<\infty.\]
Equipped with the norm $\norm{f}_{B^s_{p,q}}=\norm{f}_p+|f|_{B^s_{p,q}}$ this becomes a Banach space. 
%$\text{Lip}(s,p)$ is the set of $L^p$ functions which satisfy $\omega(f,t)_p=\mathcal{O}(t^s)$. For $0<s<1$, $\text{Lip}(s,p)=B^s_{p,\infty}$ while $\text{Lip}(1,p)=W^1_{p}$.
For more on these function spaces we refer the reader to \cite[Ch.\,5]{Stein:1970yt}.

\section{Norm Estimates}
Our approach hinges on studying the Bessel--Riesz quotient, $E_{\alpha,\mu}$, which, for $\mu>0$, defined by the multiplier
\begin{equation}\label{eqn:mmu}
m_{\alpha,\mu}(\xi):=\dfrac{\abs{\xi}^\alpha}{(\mu^2+\abs{\xi}^2)^{\frac{\alpha}{2}}}.
 \end{equation} 
We focus mainly on the case $0<\alpha\leq1$. At least two observations point to the connection with approximation theory. The first is the trivial fact that $m_{\alpha,\mu}(\xi)\to0$ pointwise as $\mu\to\infty$. The second observation starts with a formula from \cite[\S5.3.2]{Stein:1970yt}
\[\dfrac{\abs{\xi}^\alpha}{(\mu^2+\abs{\xi}^2)^{\frac{\alpha}{2}}}=\paren{1-\frac{\mu^2}{\mu^2+\abs{\xi}^2}}^{\frac\alpha2}=1-\sum_{j=1}^\infty a_{\alpha,j}(1+|\xi\mu^{-1}|^2)^{-j},\] 
for some positive coefficients with $\sum a_{\alpha,j}=1$. By Fourier inversion we obtain
\begin{equation}\label{eqn:emu}
E_{\alpha,\mu} f(x) =f(x) - T_{\alpha,\mu} f(x),
\end{equation}
where $T_{\alpha,\mu}$ has the convolution kernel $A_{\alpha,\mu}(z)$ defined as
\begin{equation}\label{eqn:amu_series}
A_{\alpha,\mu}(z):=\mu^{d}\sum_{j=1}^{\infty}a_{\alpha,j}G_{2j}(\mu z).
\end{equation}
Here $G_s(z)$ are the standard Bessel kernels. Their well known properties imply that $A_{\alpha,\mu}(z)$ is a positive, radial, integrable function with $L^1$ norm $\norm{A_{\alpha,\mu}}_1=\norm{\sum_{j=1}^{\infty}a_{\alpha,j}G_{2j}}_1=1$. 
%Note that $A_{\alpha,\mu}(z)$ is also singular at the origin since $\lim_{|z|\to0}G_{2j}(z)=\infty$ when $0<2j\leq d$. 
Evidently, $T_{\alpha,\mu}$ is an approximate identity and, by \eqref{eqn:emu}, $E_{\alpha,\mu}$ is its approximation error. 
Thus
\[
E_{\alpha,\mu} f(x)=f(x)-\int_{\R^d} A_{\alpha,\mu}(x-y)f(y)\,dy=\int_{\R^d} A_{\alpha,\mu}(x-y)(f(x)-f(y))\,dy.
\]
Minkowski's inequality and a change of variables readily show that the \emph{order of approximation}, $\norm{E_{\alpha,\mu} f}_p$, satisfies
\begin{equation*}\label{eqn:start_lp_int_est}
\norm{E_{\alpha,\mu} f}_p\leq \int_{\R^d} A_{\alpha,1}(y)\norm{f(\cdot)-f(\cdot-y/\mu)}_p\,dy\leq \int_{\R^d} A_{\alpha,1}(y)\omega(f,|y|/\mu)_p\,dy.
\end{equation*}
Since $\omega(f,t)_p\leq2||f||_p$, it follows that
\begin{equation}\label{eqn:basic_Lp_bound}
 \norm{E_{\alpha,\mu} f}_p\leq2\norm{f}_p,
 \end{equation}
which implies that $ \norm{G_{\alpha,\mu} f}_p\leq2\norm{I_\alpha f}_p$, but we improve on this bound below.

To simplify the notation we set $A_\alpha(z):=A_{\alpha,1}(z)=\sum_{j=1}^{\infty}a_{\alpha,j}G_{2j}(z)$. Some properties of $a_{\alpha,j}$ and $G_{2j}(z)$ are gathered next.

\begin{lemma}\label{lem:bessel_properties}
\leavevmode
\begin{enumerate}[(a)]
\item $G_{2j}(z)$ is positive, radial and decreasing with $\norm{G_{2j}}_{L^1}=1$. Moreover,
 $ G_{2j}(z)= \frac{1}{2^{\frac{d+2j-2}{2}}\pi^{\frac{d}{2}}\Gamma\paren{j}}K_{\frac{d-2j}{2}}(|z|)|z|^{\frac{2j-d}{2}}$, where $K_\nu(t)$ is the modified Bessel function of the third kind.
\item We have $a_{\alpha,j}>0$ and $a_{\alpha,j}\leq C_\alpha j^{-1-\frac\alpha2}$ with $\sum_{j=1}^{\infty}a_{\alpha,j} =1$.
%and $\sum_{j=1}^{\infty}(-1)^ja_j =1-\sqrt{2}$ .
\item $\int_{\R^d} G_{2j}(y)|y|^s\,dy=C_{d,s}\frac{\Gamma\paren{j+\frac{s}{2}}}{\Gamma(j)}$ for some constant $C_{d,s}$. In addition, for $0\leq s<\alpha$, $\int_{\R^d} A_\alpha(y)|y|^s\,dy$ converges.
\end{enumerate}
\end{lemma}

\begin{proof}
\hfill
\begin{enumerate}[(a)]
\item These are proved in {\sc Aronszajn--Smith} \cite[pgs. 413--421]{Aronszajn:1961nd}.
\item We use the binomial expansion $\paren{1-t}^{\alpha/2}=1-\sum_{j=1}^\infty a_{\alpha,j}t^j$, where
\[ 
a_{\alpha,j}=\abs{\binom{\alpha/2}{j}}=\frac{1}{\Gamma(-\frac\alpha2)j^{1+\frac\alpha2}}\paren{1+o(1)}\leq C_\alpha j^{-1-\frac\alpha2}.
\] 
It follows that $\sum_{j=1}^{\infty}a_{\alpha,j} t^j$ converges absolutely for $\abs{t}\leq1$. Evaluating $\paren{1-t}^{\alpha/2}$ at $t=1$ shows that $\sum_{j=1}^{\infty}a_{\alpha,j} =1$.

\item For $|\text{Re}(\nu)|<\text{Re}(\beta)$, we use the formula \cite[Eq. 10.43.19]{NIST:DLMF}:
\begin{equation}\label{eqn:identity}
 \int_0^\infty t^{\beta-1}K_\nu(t)\,dt= 2^{\beta-2}\Gamma\paren{\frac{\beta+\nu}{2}}\Gamma\paren{\frac{\beta-\nu}{2}}.
\end{equation}
A switch to spherical coordinates combined with part (a) and \eqref{eqn:identity} gives
\[\int_{\R^d} G_{2j}(y)|y|^s\,dy=\frac{2^{2-j-\frac{d}{2}}}{\Gamma(\frac{d}{2})\Gamma(j)}\int_0^\infty t^{j+\frac{d}{2}+s-1}K_{j-\frac{d}{2}}(t)\,dt=\frac{2^s\Gamma(\frac{d+s}{2})}{\Gamma(\frac{d}{2})}\cdot\frac{\Gamma\paren{j+\frac{s}{2}}}{\Gamma(j)}.\]
Since $\Gamma(x+a)\sim\Gamma(x)x^a$ for large $x$, and $a_j\leq Cj^{-1-\alpha/2}$ we have
\begin{align*}
\int_{\R^d} A_\alpha(y)|y|^s\,dy=\sum_{j=1}^\infty a_{\alpha,j}\int_{\R^d} G_{2j}(y)|y|^s\,dy&=C_{s,d}\sum_{j=1}^\infty a_{\alpha,j}\frac{\Gamma\paren{j+\frac{s}{2}}}{\Gamma(j)}\\&\leq C \sum_{j=1}^\infty j^{-\frac{(2+\alpha-s)}{2}},
\end{align*}
which converges when $0\leq s<\alpha$.
\end{enumerate}
\end{proof}

Parts (b) and (c) of Theorem \ref{thm:bessel-riesz-modulus} are essentially contained in the next result.
\leavevmode
\begin{theorem}\label{thm:weak_jackson}
Assume $1\leq p<\infty$.
\begin{enumerate}[(i)]
\item If $\alpha=1$, there is a $C>0$ depending only on $d$ such that
\[\norm{E_{1,\mu} f}_{p}\leq C\omega(f,1/\mu)_p\paren{3+2\ln\paren{\frac{||f||_p}{2\omega(f,1/\mu)_p}}}.\]
\item If $0<\alpha<1$, there is a constant $C$ depending only on $d$ and $\alpha$ such that
\[
\norm{E_{\alpha,\mu} f}_{p}\leq C\paren{\omega(f,1/\mu)_p+(\omega(f,1/\mu)_p)^\alpha\cdot||f||_p^{1-\alpha}}
\]
\end{enumerate}
\end{theorem}
\begin{proof}
Recall that $\norm{E_{\alpha,\mu} f}_p\leq \int A_{\alpha}(y)\omega(f,|y|/\mu)_p\,dy$. Let $R$ be a positive number to be chosen shortly. Since $A_{\alpha}(z)=\sum_{j=1}^{\infty}a_{\alpha,j}G_{2j}( z)$, it follows that
\begin{align*}
\norm{E_{\alpha,\mu} f}_{p}&\leq \sum_{j=1}^{\infty}a_{\alpha,j}\int_{\R^d}\omega\paren{f,|y|/\mu}_pG_{2j}(y)\,dy\\
&\leq \sum_{j\leq R}a_{\alpha,j}\int_{\R^d}\omega\paren{f,|y|/\mu}_pG_{2j}(y)\,dy+ \sum_{j>R}a_{\alpha,j}\int_{\R^d}\omega\paren{f,|y|/\mu}_pG_{2j}(y)\,dy.
\end{align*}
For the first sum we use the inequality $\omega(f,\gamma t)_p\leq(1+|\gamma|)\omega(f,t)_p$. In the second sum we use $\omega(f, t)_p\leq2||f||_p$. Altogether 
\[
\norm{E_{\alpha,\mu} f}_{p}\leq  \omega\paren{f,1/\mu}_p\sum_{j\leq R}a_{\alpha,j}\int_{ \R^d}(1+|y|)G_{2j}(y)\,dy+2\norm{f}_p\sum_{j>R}a_{\alpha,j}\int_{\R^d}G_{2j}(y)\,dy.
\]
By Lemma \ref{lem:bessel_properties}(c),
\[\norm{E_{\alpha,\mu} f}_{p}\leq c_{d}\omega\paren{f,1/\mu}_p\sum_{j\leq R}a_{\alpha,j}(1+j^{\frac12}) + 2\norm{f}_p\sum_{j\geq R}a_{\alpha,j}.\]
We can now split the argument into the two cases.
\begin{enumerate}[(i)]
\item \underline{The case $\alpha=1$}: We know that $a_{1,j}\leq cj^{-3/2}$ from Lemma \ref{lem:bessel_properties}(b) and can compare sums to integrals to deduce 
\begin{equation}\label{eqn:b4min}
\norm{E_{1,\mu} f}_{p}\leq c_d\paren{\omega\paren{f,1/\mu}_p(1+\ln R) + 2\norm{f}_pR^{-\frac{1}{2}}}.
\end{equation}
The choice $R=\paren{||f||_p\slash 2\omega(f,1/\mu)_p}^2$ minimizes \eqref{eqn:b4min} and completes the proof in this case.
\item \underline{The case $0<\alpha<1$}: Here $a_{\alpha,j}\leq c_\alpha j^{-1-\frac\alpha2}$ and this time the integral test yields
\begin{equation}
\norm{E_{\alpha,\mu} f}_{p}\leq c_{\alpha,d}\paren{\omega\paren{f,1/\mu}_p(1+R^{\frac12-\frac\alpha2}) + \norm{f}_pR^{-\frac{\alpha}{2}}}.
\end{equation}
This is minimized by $R=\paren{\dfrac{\alpha||f||_p}{(1-\alpha)\omega(f,1/\mu)_p}}^2$.
\end{enumerate}
\end{proof}

Theorem \ref{thm:bessel-riesz-modulus} (a) follows directly from the next result on the equivalence between the order of approximation  and the modulus of continuity.

\begin{theorem}\label{thm:main_equiv}
Suppose $f\in L^p(\R^d)$, $\alpha=1$ and $1<p<\infty$. Then
\begin{equation}\label{eqn:emu_1p}
\norm{E_{1,\mu} f}_p\approx\omega(f,1/\mu)_p.
\end{equation}
\end{theorem}
It turns out that \eqref{eqn:emu_1p} is implicit in {\sc Colzani} \cite{Colzani:1987pr} and {\sc Liu--Lu} \cite{Liu:1991oz}, but we give an independent proof. We first establish the equivalence between the order of approximation and a certain $K$--functional. Known relationships between $K$--functionals and the modulus of continuity allow us to complete the proof.

Following {\sc Ditzian--Ivanov} \cite{Ditzian:1993wq}, we introduce the $K$--functional
\begin{equation}\label{eqn:k_func}
K(t,f,|D|)_p:=\inf_{\substack{g\in L^p\\|D|g\in L^p}}\paren{\norm{f-g}_p+t\norm{|D|g}_p}.
\end{equation}
\begin{lemma}\label{lem:interim_equiv}
$K(1/\mu,f,|D|)_p\approx\norm{E_{1,\mu} f}_p$, for $1<p<\infty$.
\end{lemma}
\begin{proof}
First assume that both $g,|D|g\in L^p$. For $t>0$ define the ``dilated" Bessel kernel by $\J_s(x,t)=t^dG_s(tx)$. Then
\[
\norm{E_{1,\mu} g}_p=\mu^{-1}\norm{\J_1(\cdot,\mu)\star |D|g}_p\leq\mu^{-1}\norm{|D|g}_p.
\]
This combined with \eqref{eqn:basic_Lp_bound}  implies
\[\norm{E_{1,\mu} f}_p\leq\norm{E_{1,\mu}(f-g)}_p+\norm{E_{1,\mu} g}_p\leq\norm{f-g}_p+\mu^{-1}\norm{|D|g}_p.\]
Taking the infimum over such $g$ gives $\norm{E_{1,\mu} f}_p\leq K(1/\mu,f,|D|)_p$  which is one direction of the result.

We turn to the opposite inequality. Set $g=T_{1,\mu} f$. 
We will show that
 \[\mu^{-1}\norm{|D|g}_p:=\mu^{-1}\norm{|D|T_{1,\mu} f}_p\leq C\norm{f-T_{1,\mu} f}_p.\] On the Fourier transform side
\begin{align*}
\mu^{-1}\widehat{|D|T_{1,\mu} f}(\xi)=\frac{|\xi|}{\mu}\paren{1-\frac{|\xi|}{(\mu^2+|\xi|^2)^{\frac{1}{2}}}}\widehat{f}(\xi)&=\frac{\mu}{((\mu^2+|\xi|^2)^{\frac{1}{2}}+|\xi|)}\cdot\frac{|\xi|\widehat{f}(\xi)}{(\mu^2+|\xi|^2)^{\frac{1}{2}}}\\
&:=r(\xi)\widehat{E_{1,\mu} f}(\xi),
\end{align*}
and we only need show that $r(\xi)$ defines a bounded operator on $L^p$.
A direct computation shows that for $\xi\neq0$
\[\abs{\frac{\partial r}{\partial\xi_k}}=\abs{-\mu((\mu^2+|\xi|^2)^{\frac{1}{2}}+|\xi|)^{-2}\cdot\paren{\frac{\xi_k}{|\xi|}+\frac{\xi_k}{(\mu^2+|\xi|^2)^{\frac{1}{2}}}}}\leq\frac{2}{|\xi|}.\]

For any multi-index $\gamma$, this can be extended to $|\partial^\gamma r(\xi)|\leq C_\gamma|\xi|^{-|\gamma|}$.
The multiplier theorem shows that $r(D)$ is $L^p$ bounded for $1<p<\infty$. Hence, for $1<p<\infty$, we obtain $\mu^{-1}\norm{|D|T_{1,\mu} f}_p\leq C\norm{f-T_{1,\mu}f}_p$. To conclude, note that 
\[K(1/\mu,f,|D|)_p\leq\norm{f-T_{1,\mu}f}_p +\mu^{-1}\norm{|D|T_{1,\mu}f}_p\leq C\norm{E_{1,\mu} f}_p.\]
\end{proof}
We need a variant of Calderon's theorem.
\begin{lemma}\label{lem:D_|D|}
For $1<p<\infty$, $g\in W^1_p$ if and only if $g$ and $|D|g$ are in $L^p$.
\end{lemma}
\begin{proof}
Using the Riesz transforms $R_j$, we can write $D_{j}g=R_j(|D|g)$.
The boundedness of $R_j$ implies that $D_jg\in L^p$ whenever $|D|g\in L^p$ if $1<p<\infty$. Since $R_j$ is unbounded in $L^1$ and $L^\infty$, we cannot include the case $p=1$ or $\infty$.

For the converse, suppose that $g\in W^1_p$. Then $g=G_{1,1}\star h$ for some $h\in L^p$ (see \cite[\S5.3.4]{Stein:1970yt}). By definition, $\widehat{|D|g}(\xi)=|\xi|(|\xi|^2+1)^{-\frac{1}{2}}\widehat{h}(\xi)$, so that $|D|g=E_{1,1}h$. As $E_{1,1}$ is $L^p$ bounded,  $\norm{|D|g}_p<\infty$, completing the proof.
\end{proof}
We are in a position to make short work of Theorem \ref{thm:main_equiv}.
\begin{proof}[Proof of Theorem \ref{thm:main_equiv}]
We apply the result of {\sc Johnen--Scherer}, \cite{Johnen:1977gb}, on the equivalence of moduli of continuity and $K$--functionals. If we define
\[K(t, f, L^p,W^1_p)=\inf_{g\in W^1_p}\paren{||f-g||_p+t\sup_{|\gamma|=1}||D^\gamma g||_p},\]
their result is that $K(t, f, L^p,W^1_p)\approx\omega(f,t)_p$ for $1\leq p\leq\infty$. However, Lemma \ref{lem:D_|D|} shows that when $g\in W^1_p$, we have $\sup_{|\gamma|=1}\norm{D^\gamma g}_p\approx \norm{|D|g}_p$ for $1<p<\infty$.
This implies $K(t,f,|D|)\approx K(t, f, L^p,W^1_p)$, and we have shown
\[\norm{E_{1,\mu} f}_p\approx K(1/\mu,f,|D|)\approx K(1/\mu, f, L^p,W^1_p)\approx\omega(f,1/\mu)_p.\]
\end{proof}

Theorem \ref{thm:bessel-riesz-modulus} follows from Theorems \ref{thm:weak_jackson} and \ref{thm:main_equiv} and the identity $G_{\alpha,\mu}f=E_{\alpha,\mu} I_\alpha f$. It is worth pointing out that Theorem \ref{thm:bessel-riesz-modulus} (a) gives yet another characterization of Besov spaces.
\begin{corollary}\label{lem:besov}
 Fix $0<s<1$. For $1<p<\infty$ and $0<q\leq\infty$, we have 
\[|f|^q_{{B}^s_{p,q}}\approx \int_1^\infty (\mu^{s}||E_{1,\mu} f||_p)^q\,\frac{d\mu}{\mu}.\]
\end{corollary}

We end this section with Hardy space estimates for $E_{1,\mu}$. The real Hardy space, $H^p(\R^d)$, for $p>1$ coincide with the $L^p$ spaces. For $0<p\leq1$, it is a normed space of distributions. We denote the norm by $\norm{\cdot}_{H^p}$. For $f$ in $H^p(\R^d)$, $\omega(f,t)_{{H}^p}:=\sup_{|h|\leq t}||f(\cdot+h)-f(\cdot)||_{H^p}$ is the $H^p$ modulus of continuity. A thorough exposition can be found in \cite[Ch.\,3]{Stein:1993aq}. We have stated the minimum required to formulate the next result.

\begin{theorem}\label{thm:jackson_hardy}
 For $0<p\leq1$, we have $ \norm{E_{1,\mu} f}_{{H}^p}\leq C_p\omega(f,1/\mu)_{{H}^p}$.
\end{theorem}

The proof closely follows that given by {\sc Colzani} \cite{Colzani:1987pr} for the approximation error of Bochner--Riesz means. It uses a multiplier theorem for $H^p$ (see \cite[\S7.4.9]{Stein:1970yt}) along with a result on approximating $H^p$ functions by entire functions of exponential type (see \cite[Theorem 4.1]{Colzani:1987pr}). We omit the details since they do not differ substantially from \cite{Colzani:1987pr}. They can be found in an earlier version of this paper.

\section{Pointwise Estimates}
The proof of Theorem \ref{thm:bessel-riesz-local} is based on pointwise estimates of the kernel of the Bessel--Riesz quotient. Note that if $b(\xi)$ is either $m_{\alpha,\mu}(\xi)$ or $1-m_{\alpha,\mu}(\xi)$, then $b(\xi)$ satisfies 
\begin{equation}\label{eqn:mult_est}
|\partial^\beta_\xi b(\xi)|\leq C_\beta \abs{\xi}^{-|\beta|};\quad \xi\neq0.
\end{equation}
This symbol estimate already implies the $L^p$ boundedness of $E_{\alpha,\mu}$ for $1<p<\infty$, but does not give fast enough decay at infinity for the kernel. However, a more detailed analysis actually shows that for $\xi\neq0$
\begin{equation}\label{eqn:mmu_est}
|\partial^\beta_\xi b(\xi)|\leq C_{\alpha,\beta}|\xi|^{\alpha-|\beta|}(\mu^2+|\xi|^2)^{-\frac{\alpha}{2}}.
\end{equation}
This small refinement is the main ingredient in the following result.
\begin{lemma}\label{lem:kernel_estimate}
Suppose $b(\xi)\in L^\infty$ and satisfies \eqref{eqn:mmu_est} for $0<\alpha\leq1$. Then its kernel $B(x)$ satisfies
\[|\partial_x^\gamma B(x)| \leq  C_{\alpha,\gamma,d}
\begin{cases}
|2\mu x|^{-\frac\alpha2}|x|^{-|\gamma|-d}; &|\mu x|>1,\\
(|\mu x|^2+1)^{-\frac\alpha2}|x|^{-|\gamma|-d};&|\mu x|\leq1.
\end{cases}
\]
\end{lemma}
Near the origin, this is the standard estimate for Calderon--Zygmund kernels. The extra decay at infinity leads to a quantitative localization principle and Theorem \ref{thm:bessel-riesz-local}. Let us show this before proving the Lemma.

\begin{proof}[Proof of Theorem \ref{thm:bessel-riesz-local}]
We first show that if $f$ in $L^p$ vanishes near $x_0$, then $E_{\alpha,\mu} f(x_0) =\mathcal{O}(\mu^{-\frac\alpha2})$.  
By translation invariance, we may assume that $x_0=0$, and $\delta>0$ is such that $f=0$ for $|x|<\delta$. If $\mu\delta>1$, Lemma \ref{lem:kernel_estimate} applied to $m_{\alpha,\mu}(\xi)$ gives
\[
|E_{\alpha,\mu} f(0)| = \abs{\int_{|y|>\delta}K_{\alpha,\mu}(-y)f(y)\,dy}\leq \frac{C}{\mu^{\frac\alpha2}}\int_{|y|>\delta}\frac{|f(y)|}{|y|^{d+\frac{\alpha}{2}}}\,dy,
\]
where $K_{\alpha,\mu}$ is its kernel.  By H\"older's inequality,
\[|E_{\alpha,\mu} f(0)| \leq C_{d,p}\mu^{-\frac\alpha2}\delta^{-(\frac{d}{p}+\frac{\alpha}{2})}\norm{f}_p\leq C_{d,\delta,p}\mu^{-\frac\alpha2}\norm{f}_p.\]
The proof is completed by appealing to the identity $G_{\alpha,\mu}f=E_{\alpha,\mu} I_\alpha f$. 
\end{proof}

We turn now to the proof of Lemma \ref{lem:kernel_estimate}. It is based on a standard Littlewood--Paley type argument as in {\sc Stein} \cite[pgs.~241-247]{Stein:1993aq}. We include it to show the effect of replacing the H\"ormander--Mikhlin condition \eqref{eqn:mult_est} with \eqref{eqn:mmu_est}. 
\begin{proof}[Proof of Lemma \ref{lem:kernel_estimate}]
Let $1=\sum_{j\in\Z}\phi(2^{-j}\xi)$ be a Littlewood--Paley partition of unity.  Put 
\[B_{j}(x) = \int e^{ix\cdot\xi}\phi(2^{-j}\xi)b(\xi)\,d\xi\]
For any multi-indices $\beta$ and $\gamma$ we see
\[\abs{x^\beta(-i\partial_x)^\gamma B_{j}(x)}=\abs{\int x^\beta e^{ix\xi}\xi^\gamma\phi(2^{-j}\xi)b(\xi)\,d\xi}
\leq \int \abs{\partial_\xi^\beta(\xi^\gamma\phi(2^{-j}\xi)b(\xi))}\,d\xi.
\]
The product rule, \eqref{eqn:mmu_est}, and direct integration gives
\[
\abs{x^\beta(-i\partial_x)^\gamma B_{j}(x)}\leq C_{\gamma,\beta,d}2^{j(d+|\gamma|-|\beta|)}\frac{2^{\alpha(j-1)}}{(2^{2(j-1)}+\mu^2)^{\alpha/2}},
\]
which can be rearranged to the derivative estimate
\begin{equation}\label{eqn:m_est}
|\partial_x^\gamma B_{j}(x)|\leq C_{\gamma,M,d}2^{j(d+|\gamma|-M)}\frac{2^{\alpha(j-1)}}{(2^{2(j-1)}+\mu^2)^{\alpha/2}}|x|^{-M}.
\end{equation}
We now split the sum as
 \[\partial_x^\gamma B(x)=\sum_{2^{j-1}\leq|x|^{-1}}\partial_x^\gamma B_{j}(x)+\sum_{2^{j-1}>|x|^{-1}}\partial_x^\gamma B_{j}(x).\]
To estimate the first sum, set $M=0$ in \eqref{eqn:m_est} to find that
\begin{equation}\label{eqn:small_j_series}
\sum_{2^{j-1}\leq|x|^{-1}}|\partial_x^\gamma B_{j}(x)| \leq  C_{\gamma,d}\sum_{2^{j-1}\leq|x|^{-1}}\frac{2^{j(d+|\gamma|)}}{(1+\paren{\mu/2^{j-1}}^2)^{\frac\alpha2}}.
\end{equation}
When $2^{j-1}\leq|x|^{-1}$, we see that $(1+(\mu/2^{j-1})^2)^{-\frac\alpha2}\leq (1+|\mu x|^2)^{-\frac\alpha2}$ and summing the geometric series \eqref{eqn:small_j_series} we obtain
\begin{equation}\label{eqn:small_j}
\sum_{2^{j-1}\leq|x|^{-1}}|\partial_x^\gamma B_{j}(x)|\leq \frac{C_{\gamma,d}}{(1+|\mu x|^2)^{\frac\alpha2}}\frac{1}{|x|^{d+|\gamma|}}.
\end{equation}

For the second sum, we set $M$ to be the smallest integer greater than $|\gamma|+d+1/2$ and arrive at
\begin{align}\label{eqn:big_j_series}
\sum_{2^{j-1}>|x|^{-1}}|\partial_x^\gamma B_{j}(x)| &\leq  C_{\gamma,d,M}|x|^{-M}\sum_{2^{j-1}>|x|^{-1}}2^{j(d+|\gamma|-M)}\frac{2^{\alpha(j-1)}}{\paren{\mu2^{j-1}\paren{\frac{\mu}{2^{j-1}}+\frac{2^{j-1}}{\mu}}}^{\frac\alpha2}}\nonumber\\
&\leq C_{\gamma,d,M}|x|^{-M}\mu^{-\alpha/2}\sum_{2^{j-1}>|x|^{-1}}\frac{2^{j(d+|\gamma|+\frac{\alpha}{2}-M)}}{\paren{\frac{\mu}{2^{j-1}}+\frac{2^{j-1}}{\mu}}^{\frac\alpha2}}.
\end{align}
Setting $t=\mu/2^{j-1}$ and $L=|\mu x|$. If $2^{j-1}>|x|^{-1}$, we see that $0<t\leq L$. A direct calculation shows
\[\sup\limits_{0<t\leq L}(t+t^{-1})^{-\frac\alpha2}=
\begin{cases}
2^{-\frac\alpha2}; &L>1,\\
(L+L^{-1})^{-\frac\alpha2};&L\leq1.
\end{cases}
\]
We use this to sum the geometric series in \eqref{eqn:big_j_series} and obtain 
\begin{equation}\label{eqn:big_j}
\sum_{2^{j-1}>|x|^{-1}}|\partial_x^\gamma B_{j}(x)|\leq  C
\begin{cases}
|2\mu x|^{-\frac\alpha2}|x|^{-|\gamma|-d}; &|\mu x|>1,\\
(|\mu x|^2+1)^{-\frac\alpha2}|x|^{-|\gamma|-d};&|\mu x|\leq1.
\end{cases}
\end{equation}
Combining \eqref{eqn:big_j} with the earlier estimate \eqref{eqn:small_j} completes the proof.
\end{proof}
A similar argument establishes a  H\"ormander--type condition which we include here for completeness.
\begin{corollary}\label{cor:horm_cond}
If $b(\xi)$ satisfies \eqref{eqn:mmu_est}, then its kernel satisfies
\begin{equation}
\int\limits_{|x|\geq2|y|}\abs{B(x+y)-B(x)}\,dx\leq C 
\begin{cases}
|2\mu y|^{-\alpha/2}; &|\mu y|>1,\\
(|\mu y|^2+1)^{-\alpha/2};&|\mu y|\leq1.
\end{cases}
\end{equation}
\end{corollary}
We do not know if Lemma \ref{lem:kernel_estimate} is sharp, but better kernel estimates near infinity will lead to better localization and pointwise estimates.

\subsection*{Acknowledgement}
I would like to thank Professor L. Colzani for clarifying some points in \cite{Colzani:1987pr}. Professor A. Larrain--Hubach also made helpful comments on an earlier draft. Any errors are, of course, mine. 
%\subsection*{Availability of Data and Material}
%Not applicable.
%\subsection*{Conflict of Interest}
%The author states that there is no conflict of interest.
\thispagestyle{empty}
%\newpage
\bibliography{lib_papers} %This is the global bib library
\bibliographystyle{amsplain}
\end{document}